\documentclass[journal]{new-aiaa}
\usepackage[utf8]{inputenc}

\usepackage{graphicx}
\usepackage{amsmath}
\usepackage[version=4]{mhchem}
\usepackage{siunitx}
\usepackage{longtable,tabularx}
\title{Autonomous Time-Optimal Many-Revolution Orbit Raising for Electric Propulsion GEO Satellites via Neural Networks}

\author{Haiyang Li \footnote{PhD Candidate, School of Aerospace Engineering, lihy15@mails.tsinghua.edu.cn; Visiting PhD Candidate, Department of Aerospace Science and Technology, haiyang.li@mail.polimi.it}}
\affil{Tsinghua University, Beijing, China, 100084}
\affil{Politecnico di Milano, Milan, Italy, 20156}
\author{Francesco Topputo\footnote{Associate Professor, Department of Aerospace Science and Technology, francesco.topputo@polimi.it.}}
\affil{Politecnico di Milano, Milan, Italy, 20156}
\author{Hexi Baoyin\footnote{Professor, School of Aerospace Engineering; baoyin@tsinghua.edu.cn. Senior member AIAA}}
\affil{Tsinghua University, Beijing, China, 100084}

\begin{document}

\maketitle

\section{Introduction}
\lettrine{E}{lectric} propulsion (EP) Geostationary Earth orbit (GEO) satellites are more efficient because they cost a lower fuel consumption and thus increase the payload mass or reduce the launch mass, compared to conventional chemical propulsion GEO satellites. EP has been tested and applied a lot in many deep space missions, and its application to the Earth orbit satellites has been a growing interest \cite{topputo2019catalogue, ceccherini2018system}. The world-first all-electric GEO platform, Boeing 702SP, has already been launched in March 2015 by a Falcon 9 rocket. China also launched its first EP-GEO satellite in April 2017. The design of GEO satellites is rapidly changing from classic high-thrust propulsion more and more toward low-thrust propulsion. Low-thrust orbit raising brings in a many-revolution long-duration transfer that is a special characteristic of EP-GEO satellites.

Autonomous orbit raising is necessary for future EP-GEO missions because the support from the ground is expensive when transfer time is very long, for instance, hundreds of days for low-thrust orbit raising from low Earth orbit (LEO) to GEO. 
Recently, neural networks (NNs) have attracted the interest of many researchers and have had applications in various fields in recent years, including astrodynamics \cite{izzo2018survey}.
For landing problems and short-duration orbit transfers, NNs are trained to learn the optimal state-control pairs  \cite{sanchez2018real,cheng2018real,izzo2019machine} and image-control relations \cite{furfaro2018deep}.
Their results show that NNs have excellent performance in learning the optimal control of short-duration problems. In this Note, NNs are trained to learn the optimal control of many-revolution long-duration problems, that is, time-optimal low-thrust orbit raising problems.

Optimization of time-optimal low-thrust orbit raising problems is necessary before training neural networks. 
The many-revolution long-duration characteristic makes it extremely difficult for numerical optimization methods to converge.
Some methods based on Edelbaum's approach \cite{edelbaum1961propulsion} have been proposed under the assumption of quasi-circular orbits \cite{casalino2007improved, kluever2011using}. However, Edelbaum's approach obtains only sub-optimal solutions, and the quasi-circular orbit assumption limits its application.
Besides Edelbaum's approach, some control laws are designed to implement low-thrust orbit raising tasks \cite{gao2007near, guelman2016closed}, but their optimality cannot be ensured.
Indirect methods are applied a lot in optimizing low-thrust trajectories \cite{zhang2015low, jiang2017improving, chen2018multi, yang2018fast} because of its optimality. In indirect methods, the original optimal control problem is transformed into a two-point boundary value problem (TPBVP) which can be solved by shooting methods. It is difficult for the shooting of low-thrust orbit raising problems to converge because of the many-revolution long-duration characteristic. Effective techniques such as homotopy, switching function detection, and computing the analytic state transition matrix are applied to improve the convergence rate \cite{jiang2012practical, chi2018power, gonzalo2017unified}. In addition, modified equinoctial elements (MEE) are utilized in \cite{ gonzalo2017indirect, li2017j2} to overcome difficulties in convergence.

In this Note, a number of trajectories are optimized to generate datasets that are used for NNs to learn. The trajectory starting from a nominal departure orbit is optimized at first using the method proposed in \cite{topputo2019catalogue}, and then departure orbits are randomly created around the nominal departure orbit and corresponding trajectories are optimized using the solutions of the nominal departure orbit as the initial guess.
As NN performance is highly affected by hyper-parameters, random search \cite{bergstra2012random} is applied to search for hyper-parameters because random search is effective and easy to implement \cite{li2019approximators}.
As will be shown, NN-driven control is capable for time-optimal low-thrust autonomous orbit raising. This work shows that besides the short-duration optimal control problems such as landing problems that are studied mostly, NNs can also be applied in this very long duration optimal control problem.

The rest of this Note is organized as follows. In Section II, trajectory optimization methods and dataset generation are presented. In Section III, the architecture of NNs and random search is detailed, and how to evaluate NN performance is also defined. Simulation results are shown in Section IV and conclusions are drawn in Section V.

\section{Trajectory Optimization}

\subsection{Time-Optimal Low-Thrust Problem}
Time-optimal low-thrust trajectories are optimized with indirect methods. The orbital states are expressed in modified equinoctial elements (MEE; $p$, $e_x$, $e_y$, $h_x$, $h_y$, $L$). The relations between MEE and the classical orbital elements are

\begin{equation}
  \label{coe2mee}
    \begin{split}
    &p=a(1-e^2)\\
    &e_x=e\cos(\omega+\Omega)\\
    &e_y=e\sin(\omega+\Omega)\\
    &h_x=\tan(i/2)\cos(\Omega)\\
    &h_y=\tan(i/2)\sin(\Omega)\\
    &L=\omega+\Omega+\theta
    \end{split}
\end{equation}
where $a$ is the semi-major axis, $e$ the eccentricity, $i$ the orbital inclination, $\Omega$ the right ascension of the ascending node, $\omega$ the perigee anomaly, $\theta$ the true anomaly; $p$ is the semilatus rectum, and $L$ the true longitude.

The dynamic equations of low-thrust propelled spacecraft in a two-body model are:
\begin{equation}
\label{dyn}
\begin{split}
    &\dot{\textbf{\textit{x}}} = u \dfrac{T_{\max }}{m} \textbf{\textit{M}} \boldsymbol{\alpha}  + \textbf{\textit{D}},\\
    &\dot m =  - u\dfrac{T_{\max}}{I_{sp}{g_0}}
\end{split}
\end{equation}
where \textbf{\textit{x}} = [\textit{p}, $e_x$, $e_y$, $h_x$, $h_y$, \textit{L}], \textit{m} is the instantaneous mass of the spacecraft, $g_0$ the standard value of gravitational acceleration, $T_{\max}$ the maximum achievable thrust magnitude, and $I_{sp}$ the thruster specific impulse. \textbf{\textit{M}} is a transformation matrix and \textbf{\textit{D}} the gravity vector; their expression can be found in \cite{gonzalo2017indirect}. $u\in[0,1]$ is the engine thrust ratio, and $\boldsymbol{\alpha}$ is the thrust unit vector of the thrust
\begin{equation}
    \boldsymbol{\alpha}=[\cos\beta \cos\alpha, \cos\beta \sin\alpha, \sin\beta]^{\top}
\end{equation}
where $\alpha$ is the azimuth angle and $\beta$ the elevation angle. The thrust vector is defined in the Local Vertical-Local Horizontal (LVLH) frame. 

In minimum-time problems, the performance index takes the forms:

\begin{equation}
    J=\int_{t_0} ^{t_f} 1 \textrm{d} t
\end{equation}

Introducing the co-state vector $\boldsymbol{\lambda}$(\textit{t})= [$\boldsymbol{\lambda_x}$, $\lambda_m$], the Hamiltonian is
\begin{equation}
\label{H}
H = u \dfrac{T_{\max }}{m} \boldsymbol{\lambda}^\top_x \textbf{\textit{M}} \boldsymbol{\alpha} + \boldsymbol{\lambda}^\top_x \textbf{\textit{D}} - u \dfrac{T_{\max}\,  \lambda_m}{I_{sp}\,{g_0}} + 1
\end{equation}
and the dynamics of $\boldsymbol{\lambda}$ is given by
\begin{equation}
\label{dyn_costate}
\dot{\boldsymbol{\lambda}} = -\dfrac{\partial H}{\partial \boldsymbol{x}}
\end{equation}

Applying Pontryagin’s minimum principle, the direction of the optimal thrust is determined as:

\begin{equation}
    \label{alpha}
    \boldsymbol{\alpha}^*=-\frac{\boldsymbol{M}^\top\boldsymbol{\lambda}_{mee}}{\|\boldsymbol{M}^\top\boldsymbol{\lambda}_{mee}\|}
\end{equation}

while the optimal thrust magnitude is
\begin{equation}
\label{utop}
\begin{split}
    &u = \left\{ {\begin{array}{*{20}{l}}
{0, \qquad {\rm{if}} \ \rho > 0},\\
{1, \qquad {\rm{if}} \ \rho < 0},\\
{[0,1], \ \ \! {\rm{if}} \ \rho = 0},
\end{array}} \right.\\
\\
    &\rho =  - \dfrac{I_{sp}g_0\|\textbf{\textit{M}}^\top \boldsymbol{\lambda_x}\|}{m} - \lambda_m
\end{split}
\end{equation}

We notice that $\lambda_m>0$ since $\dot{\lambda}_m=-\partial H/\partial m <0$ and $\lambda_m(t_f)=0$. Therefore, $\rho$ in Eq.(\ref{utop}) is always negative, and thus $u=1$ for all times in minimum-time problems. Eventually, $t_f$ is found through the transversality condition $H(t_f)=0$.

Given the boundary condition, the optimal control problem can be transformed into a TPBVP and it can be solved by shooting methods solving shooting equations. In our problem, the final mass is set to be free, and the final true longitude $L$ can be both free or fixed. 
According to the optimal control theory, if one of final states is free, then the corresponding co-state at final time has to be zero. The boundary vector $B_x$ and $B_{\lambda}$ are given as $B_{x}^{L-free}=[1,1,1,1,1,0,0]$ and $B_{\lambda}^{L-free}=[0,0,0,0,0,1,1]$ or $B_{x}^{L-fixed}=[1,1,1,1,1,1,0]$ and $B_{\lambda}^{L-fixed}=[0,0,0,0,0,0,1]$.
Then the optimal control problem can be transformed into a TPBVP and can be solved by shooting methods solving shooting equations:

\begin{equation}
    \Phi(\boldsymbol{\lambda}_0, t_f)=[(x_f-x(t_f)) \circ B_x,(\lambda_f-\lambda(t_f)) \circ B_{\lambda},H(t_f)]^{\top}=0
\end{equation}
with shooting variables $(\boldsymbol{\lambda}_0,t_f)$, where $\circ$ is the Hadamard product.

\subsection{Thrust Homotopy}
In solving the TPBVP, finding an initial guess
for $\boldsymbol{\lambda}_0$ and $t_f$ that assure convergence is difficult in this many-revolution low-thrust optimization problem, therefore a continuation method named thrust homotopy is employed \cite{topputo2019catalogue}. In thrust homotopy, a high maximum thrust is used at first, and then gradually decrease the maximum thrust to the nominal maximum thrust. For the starting high maximum thrust, random guess is used for the shooting. In the homotopy process, the relationship between $t_f$ and $T$ is  $t_f\times T \approx$ const as found in \cite{caillau2012minimum}. Thus, the guess of $t_f$ is given as this:
\begin{equation}
    t_f(T_{i+1})=T_i\times t_f ^*(T_i)/T_{i+1}
\end{equation}
where $t_f ^*(T_i)$ is the optimal transfer time when using a thrust magnitude $T_i$.

To guess $\boldsymbol{\lambda}_0 ^*$, the continuation method is selected between zero path:

\begin{equation}
    \boldsymbol{\lambda}_0(T_{i+1})=\boldsymbol{\lambda}_0 ^*(T_i)
\end{equation}
and linear path:
\begin{equation}
    \boldsymbol{\lambda}_0(T_{i+1})=\boldsymbol{\lambda}_0 ^*(T_i)+\frac{T_{i+1}-T_i}{T_i-T_{i-1}}\left(\boldsymbol{\lambda}_0 ^*(T_i)-\boldsymbol{\lambda}_0 ^*(T_{i-1})\right)
\end{equation}

The first attempt is zero path but linear path will be carried out in case that it does not converge for some situations such as low thrust-to-mass ratios and high eccentricities of the starting orbit. 

In this work, the optimal trajectories are found by using Low-Thrust Trajectory Optimizer (LT2O), a tool developing at Politecnico di Milano \cite{gonzalo2017indirect,zhang2015low}.

\subsection{Dataset Generation}
A time-optimal low-thrust trajectory from a nominal departure orbit is optimized using methods proposed above at first. Two types of departure orbits, low Earth orbit (LEO) and geostationary transfer orbit (GTO) are considered. The departure LEO is a circular orbit with a semi-major axis of 6671 km and an inclination of 5 deg. The starting GTO is a high eccentric orbit with a perigee radius of 6671 km, an apogee radius of 42164 km, and an inclination of 5 deg. The inclination is compatible with launches from Kourou. Both the $\Omega$ and $\omega$ are set to zero without any loss of generality. All transfers start from the apogee. The final state is GEO with a semi-major axis of 42164 km and inclination of 0 deg. The final true longitude is free. The initial mass of the spacecraft is 1000 kg. Two types of thruster, one with maximum thrust 1 N and Isp 2000 s and the other with maximum thrust 0.2 N and Isp 2500 s, are considered. Therefore, altogether four nominal trajectories are optimized.

To generate datasets, departure orbits are created randomly around the nominal departure orbit:

\begin{equation}
\label{starting}
\begin{split}
    &{\textbf{\emph{r}}} = {\textbf{\emph{r}}}_n+c_r\Delta r_{max}{\textbf{\emph{a}}}_r\\
    &{\textbf{\emph{v}}} = {\textbf{\emph{v}}}_n+c_v\Delta v_{max}{\textbf{\emph{a}}}_v
\end{split}
\end{equation}
where ${\textbf{\emph{r}}}_n$ and ${\textbf{\emph{v}}}_n$ are the nominal initial position and velocity, $c$ random numbers between [0 1], and ${\textbf{\emph{a}}}$ random unit vectors. Considering the orbit determination error, the max errors are set to $\Delta r_{max}$= 100 m and $\Delta v_{max}$ = 0.1 m/s.
The initial guesses of random departure orbits are the solutions of the nominal departure orbits.

In each case, 1,000 trajectories are generated as the training dataset, 100 trajectories as the validation dataset, and 100 trajectories as the test dataset. In each trajectory 1,000 points are randomly selected into the dataset. At each point, states, that are mass and MEE, and action, represented by two direction angles $\alpha$ and $\beta$, are stored.

\section{Neural Network Design}

\subsection{Architecture of NNs}

The structure of a feed-forward neural network with multiple hidden layers is determined by the number of layers $n_{layer}$ and the number of neurons $n_{neuron}$ at each hidden layer. At each layer, denoting the input \textbf{\textit{l}}$_i$, the output \textbf{\textit{l}}$_{i+1}$ is calculated as follows:

\begin{equation}
\label{layer}
\textbf{\textit{l}}_{i+1} = \textbf{\textit{f}}( \textbf{\textit{w}}{\textbf{\textit{l}}_i} + \textbf{\textit{b}})
\end{equation}
where \textbf{\textit{w}} is the weight matrix, \textbf{\textit{b}} the bias vector, and \textbf{\textit{f}} a nonlinear function named the activation function. There are three most commonly used activation functions expressed in Eq.\eqref{acfun} for the hidden layers: the sigmoid function (sig), the hyperbolic tangent (tanh) function, and the rectified linear (ReLu) function.

\begin{equation}
\label{acfun}
\begin{split}
    &f_{\rm sig}(x) = \dfrac{1}{1 + {e^{ - x}}},\\
    &f_{\rm tanh}(x) = \dfrac{e^x - e^{ - x}}{e^x + e^{ - x}},\\
    &f_{\rm relu}(x) = \max (0,x)
\end{split}
\end{equation}

The training is to adjust the value of the parameters of each layer to minimize the loss function in the form of mean squared error:

\begin{equation}
\label{loss}
E = \dfrac{1}{n}\sum\limits_{i = 1}^n \left( \hat y - y \right)^2,
\end{equation}
where y is the actual value,  $\hat y$ the estimated value, and $n$ the number of data in the training iteration.  
Gradient descent (GD) algorithms are state-of-art in training the parameters, e.g.,

\begin{equation}
\label{GD}
\textbf{\textit{w}}' = \textbf{\textit{w}} - \eta \dfrac{\partial E}{\partial \textbf{\textit{w}}},
\end{equation}
where $\eta$ is the learning rate. Some modified gradient descent algorithms are also very effective, such as momentum gradient descent (MGD) \cite{sutskever2013importance} and Adam gradient descent (AGD) \cite{kingma2014adam}.

NN performance is highly sensitive to hyper-parameters. As shown in Eq.\eqref{GD} the learning rate $\eta$ is a important factor that affect the training. With a small learning rate, the NN may take a very long time to converge, while with a large learning rate the NN may oscillate around the optimal solution. Therefore, it is a good strategy to start from a large learning rate and decrease the learning rate gradually. Two decay schedules are applied in this paper, exponential decay (ED) and natural exponential decay (NED):
\begin{equation}
\label{ED}
\begin{split}
    &\eta_{ED} ' = \eta {c ^t},\\
    &\eta_{NED} ' = \eta {e^{ - c t}},
\end{split}
\end{equation}
where \textit{c} is the decay rate and \textit{t} is the step.

In training, the initialization of weights is also a factor that can affect NN performance. Two initialization schemes are used: Fan$\_{in}$ initializer and Fan$\_{avg}$:
\begin{equation}
\begin{split}
    &x_{in}=\sqrt{\dfrac{2}{in}},\\
    &x_{avg}=\sqrt{\dfrac{6}{in+out}}
\end{split}
\end{equation}
where $in$ and $out$ are the number of units of the previous and following layers.

Besides the hyper-parameters mentioned above, the batch size B is also a hyper-parameter that must be chosen. A batch of data is put into the DNN to train and then the weights are updated. Early stopping is utilized here, and the training will stop when there is no improvement in the last N epochs.

The form of input, which is also known as features, is also an important factor that affects the performance of NNs. For the problem considered in this Note, the input is a set of transfer parameters. Three forms are considered: modified equinoctial elements MEE ($p$, $e_x$, $e_y$, $h_x$, $h_y$, $L$), classical orbit elements COE ($a$, $e$, $i$, $\Omega$, $\omega$, $f$), and Cartesian form ECI ($x$, $y$, $z$, $v_x$, $v_y$, $v_z$), and the corresponding input forms are denoted F$_{mee}$, F$_{coe}$, and F$_{eci}$, respectively.
All of the inputs are normalized to make the average of inputs is zero and standard deviation of inputs is one.

\subsection{Selection of NN models}
The NN performance is highly affected by hyper-parameters. Many model choices are implemented by manual search, which gives the impression of neural network training as an art \cite{bengio2012practical}. Random search \cite{bergstra2012random} is applied in this paper to search for hyper-parameters, as suggested in \cite{li2019approximators} because it is effective and easy to implement. Parameters considered in this Note are: number of layers ($n_{\rm layer}$), number of neurons at each hidden layer ($n_{\rm neuron}$), activation function ($f$), weights initializer (${\rm ini}$), batch size ($B$), optimizer (${\rm opt}$), initial learning rate $\eta$, decay model (${\rm dm}$), decay step (${\rm ds}$), and decay rate ($c$), and input form $F$. The search space is listed in Table \ref{tab_RS}. These parameters are uniformly random, except that the learning rate is uniform in the log-domain. The step for early stopping is set to 50. It should be noted that the minimum layer number can be one, which makes up a shallow network. To evaluate the performance of the shallow network with many neurons, when the layer number is one, the neuron number is set to be twice the value selected from the search space.

\begin{table}[htbp]
	\fontsize{10}{10}\selectfont
	\caption{Search space of random search}
	\label{tab_RS}
        \centering 
    \begin{tabular}{c | r} 
        \hline
        Hyper-parameters& Search space\\\hline
        $n_{\rm layer}$& 1--7\\
        $n_{neuron}$& 32--512\\
        $f$& sigmoid, tanh, relu\\
        ${\rm ini}$& Fan$\_{in}$, Fan$\_{avg}$\\
        $B$& 100 -- 1000\\
        ${\rm opt}$& GD, MGD, AGD\\
        $\eta$& 0.1 -- 0.0001\\
        ${\rm dm}$& ED, NED\\
        ${\rm ds}$& 100 -- 500\\
        $c$& 0.8 -- 1\\
        $F$& $F_{mee}$,$F_{coe}$,$F_{eci}$\\
        \hline
    \end{tabular}
\end{table}

The selection of NN models is based on the validation dataset. The test dataset is finally used to evaluate the performance of the selected model. When a model-selection process is considered the validation and test datasets must be separated.

\subsection{Evaluation of NN performance}
To check the accuracy of NN predictions, mean absolute percentage error (MAPE) is used on the test dataset:
\begin{equation}
\label{error}
{\rm{MAPE}} = \dfrac{1}{n}\sum\limits_{i = 1}^n {\dfrac{{\left| {\hat y - y} \right|}}{y}}
\end{equation}
where y is the actual value,  $\hat y$ the estimated value, and $n$ the number of data in the test dataset.

NN-driven trajectories are integrated using the same dynamics in Eq.\eqref{dyn} except that the thrust pointing angle is given by a NN controller rather than optimal control.
The NN-driven trajectory will reach the vicinity of GEO, and the distance  between the NN-driven trajectory and GEO is defined as the velocity increment needed from the current trajectory to GEO estimated by Edelbaum's approach \cite{edelbaum1961propulsion}:
\begin{equation}
    \Delta v = v_0\sqrt{(0.5\dfrac{\delta a}{a_0} + 0.649\delta e)^2 + (1.571\delta i)^2} 
\end{equation}
where $\delta$ is the orbital element difference between the current trajectory and GEO, $a_0$ and $v_0$ the semi major axis and velocity of the current trajectory. It should be noted that this $\Delta v$ is just a measurement and not the actual $\Delta v$ needed because actual $\Delta v$ can be lower after optimization \cite{topputo2011optimal}.

The NN-driven trajectory stops when the $\Delta v$ begins to increase. The final semimajor-axis $a_f$, eccentricity $e_f$, and inclination $i_f$ can be used to evaluate the NN performance.

To evaluate the NN performance in terms of optimality, it is fair to compare their performance indexes $J$ using the same initial and final states. 
The indirect method is applied using the same initial and final states of NN-trajectory. In this case the final true longitude $L$ is fixed. The optimality is evaluated also by MAPE of their performance indexes $J$.

\section{Simulation} 

\subsection{Results of nominal trajectories}
\begin{figure}[!htb]
\centering
\includegraphics[width=1.0\textwidth]{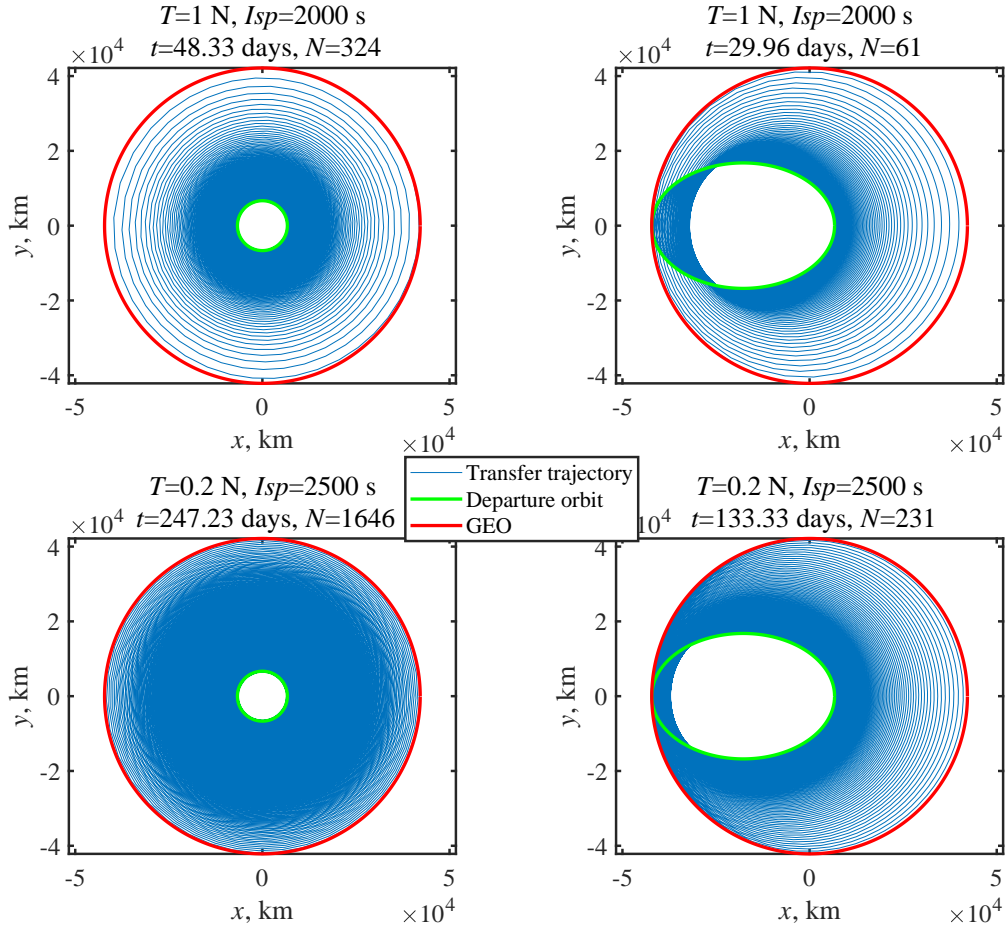}
\caption{Nominal trajectories.}
\label{fig:tra}
\end{figure}

Optimization results of nominal trajectories are presented firstly. In thrust homotopy, the starting maximum thrust is set to be 50 N, and then the maximum thrust gradually reduces to 1 N in Isp 2000 s cases and to 0.2 N in Isp 2500 s cases.
Trajectories starting from LEO or GTO are shown in Fig.\ref{fig:tra}. The transfer days $t$ and revolution number $N$ are also given in the figure. The longest duration case is the one that starts from LEO with thrust 0.2 N and Isp 2500 s, and its transfer time is 247.23 days and revolution number is 1646.

\begin{figure}[!htb]
    \centering
    \includegraphics[width=1.0\textwidth]{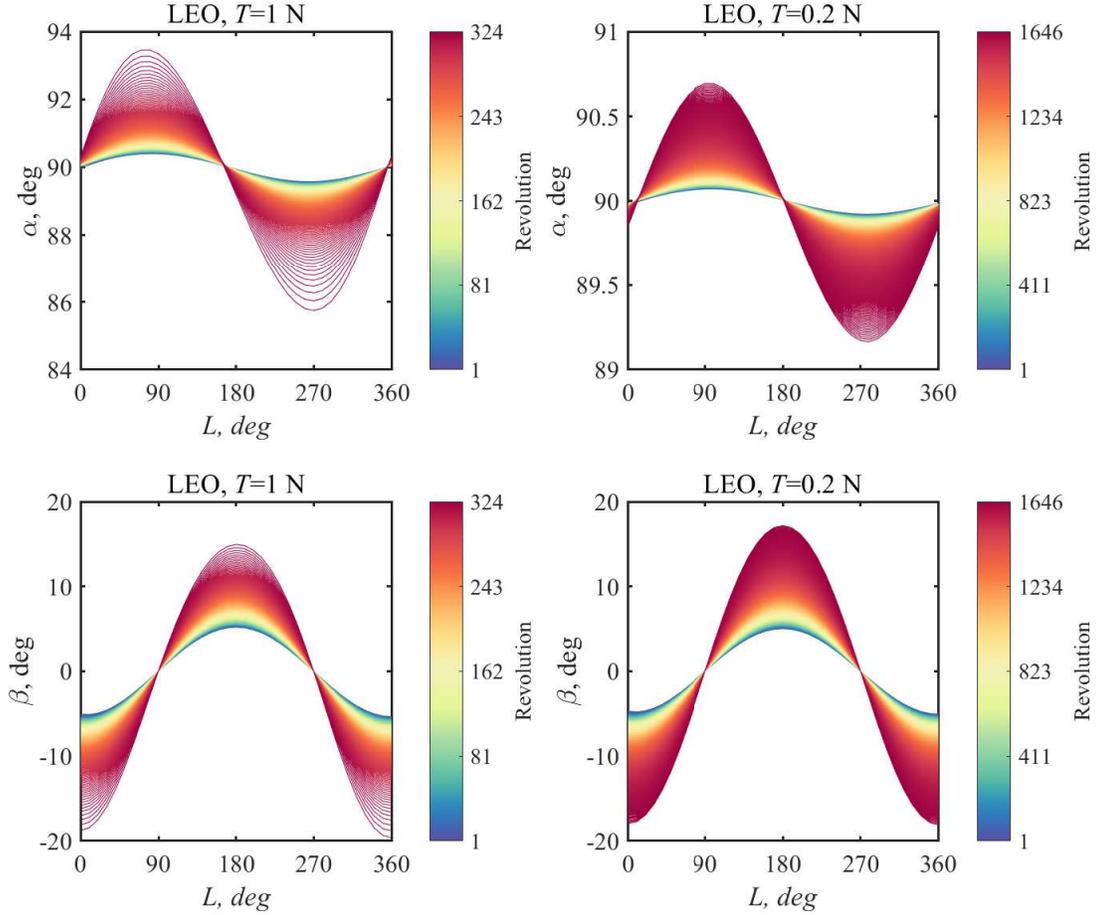}
    \caption{Control angles of LEO cases}
    \label{fig:LeoAngle}
\end{figure}

The control angles changing with true longitude $L$ for LEO cases are shown in Fig.\ref{fig:LeoAngle}, where the true longitude $L$ is modified between 0 and 360 deg and the colorbar represents the revolution number. In LEO cases, the in-plane azimuth angles $\alpha$ is centralized near 90 deg to raise the orbit altitude. The out-plane elevation angles $\beta$ distributes between about [-20, 20] deg symmetrically are used to reduce the inclination, and maximum of $\beta$ becomes higher as the revolution increases because it is more efficient to reduce the inclination at the high orbit.

\begin{figure}[!htb]
    \centering
    \includegraphics[width=1.0\textwidth]{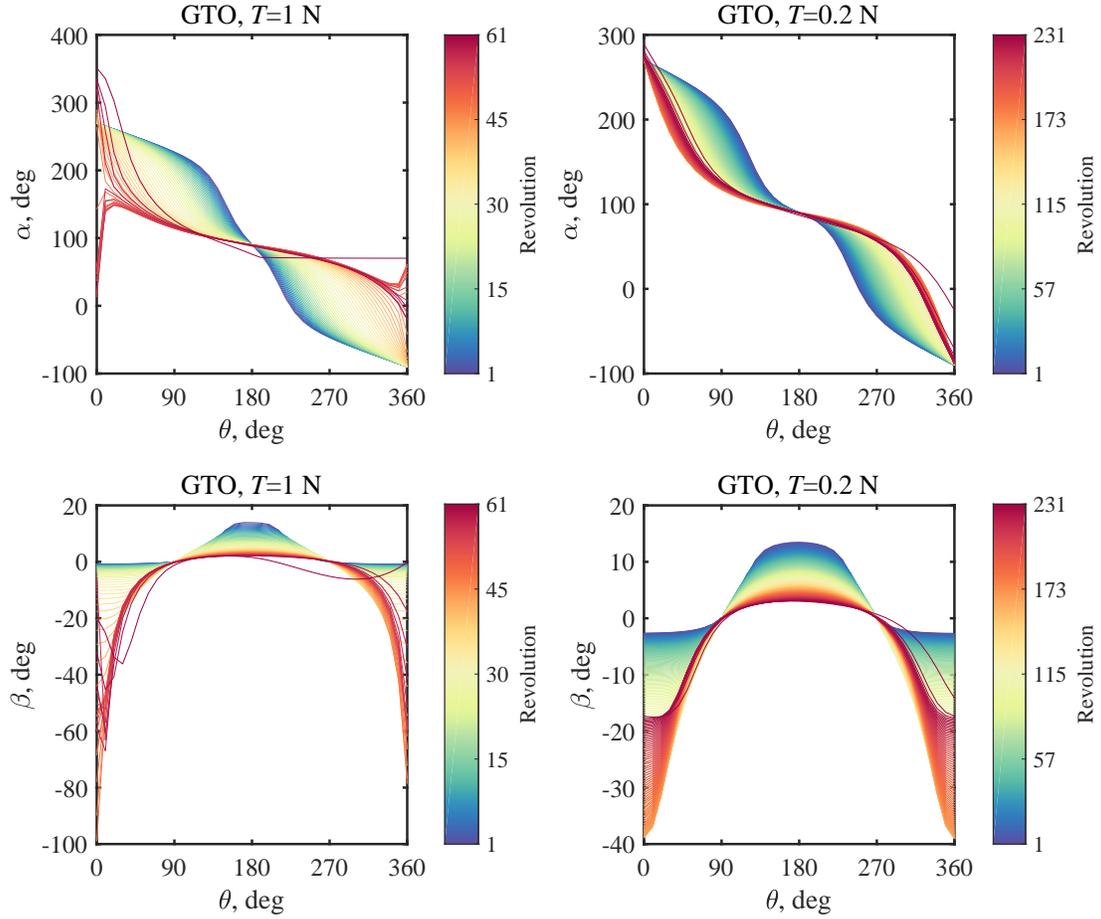}
    \caption{Control angles of GTO cases}
    \label{fig:GtoAngle}
\end{figure}

The control angles changing with true anomaly $\theta$ for GTO cases are shown in Fig.\ref{fig:GtoAngle}. In GTO cases, $\alpha$ and $\beta$ behave differently from that in LEO cases.
The in-plane azimuth angles $\alpha$ changes from 0 to 360 deg within one revolution to reduce the eccentricity and raise the altitude. $\alpha$ is about 270 deg when $\theta$ is around 0 deg (near perigee) and about 90 deg when $\theta$ is around 180 deg (near apogee), which is efficient to reduce the eccentricity, and at the same time the altitude of apogee is also reduced as shown in Fig.\ref{fig:tra}. As the revolution number increases, $\alpha$ is near 90 deg more to raise the altitude.
The out-plane elevation angles $\beta$ is no longer symmetrically distributed. In the first tens of revolutions $\beta$ is higher when $\theta$ is around 180 deg (near apogee), because it is more efficient to reduce the inclination at apogee for a high eccentric orbit. In the last tens of revolutions when the eccentricity is already reduced, $\beta$ is higher when $\theta$ is around 0 deg (near perigee), because the velocity in the perigee is higher which makes it more difficult to change the inclination, more thrust should be allocated to the out-plane part.

These results shows that the many-revolution time-optimal low-thrust transfer can be optimized effectively using the indirect method proposed in Section II.
It should be noted that although some simple laws can be concluded, the actual value of control angles should still be determined by the optimization algorithm. In the followings, NNs are trained to learn the optimal control.

\subsection{Results of neural networks}
The datasets are generated using the approach described in Section II. For each case, NN models are searched using random search respectively. The iteration of hyper-parameter random search is set to 50.
Some selected hyper-parameters and MAPE of control angles are given in Table \ref{tab_hyper}.
It can be found that the accuracy in LEO cases is higher than that in GTO cases in general because the changes of control angles in LEO cases are simple while they are quite complex in GTO cases. In addition, all selected input features are Cartesian form.
NN models with the best results of each case will be used in the following analysis.
\begin{table}[htbp]
	\fontsize{10}{10}\selectfont
	\caption{Hyper-parameters and training results of neural networks}
	\label{tab_hyper}
        \centering 
    \begin{tabular}{l | r | r | r |r |r |r |r} 
        &$F$&$n_{layers}$ & $n_{neurons}$ & $f$ & $opt$ & MAPE$\alpha$  & MAPE$\beta$ \\
        \hline
        LEO, T=1 N& ECI & 6 & 350 & tanh & AGD & 0.0113\% & 0.5453\%\\
        LEO, T=0.2 N& ECI & 6 & 452 & relu & AGD & 0.0137\% & 1.6961\%\\
        GTO, T=1 N& ECI & 4 & 176 & sigmoid & AGD & 0.4769\% & 3.8810\%\\
        GTO, T=0.2 N& ECI & 5 & 138 & sigmoid & AGD & 2.1169\% & 2.5237\%\\
        \hline
    \end{tabular}
\end{table}

\begin{figure}[h]
\centering
\includegraphics[width=1.\textwidth]{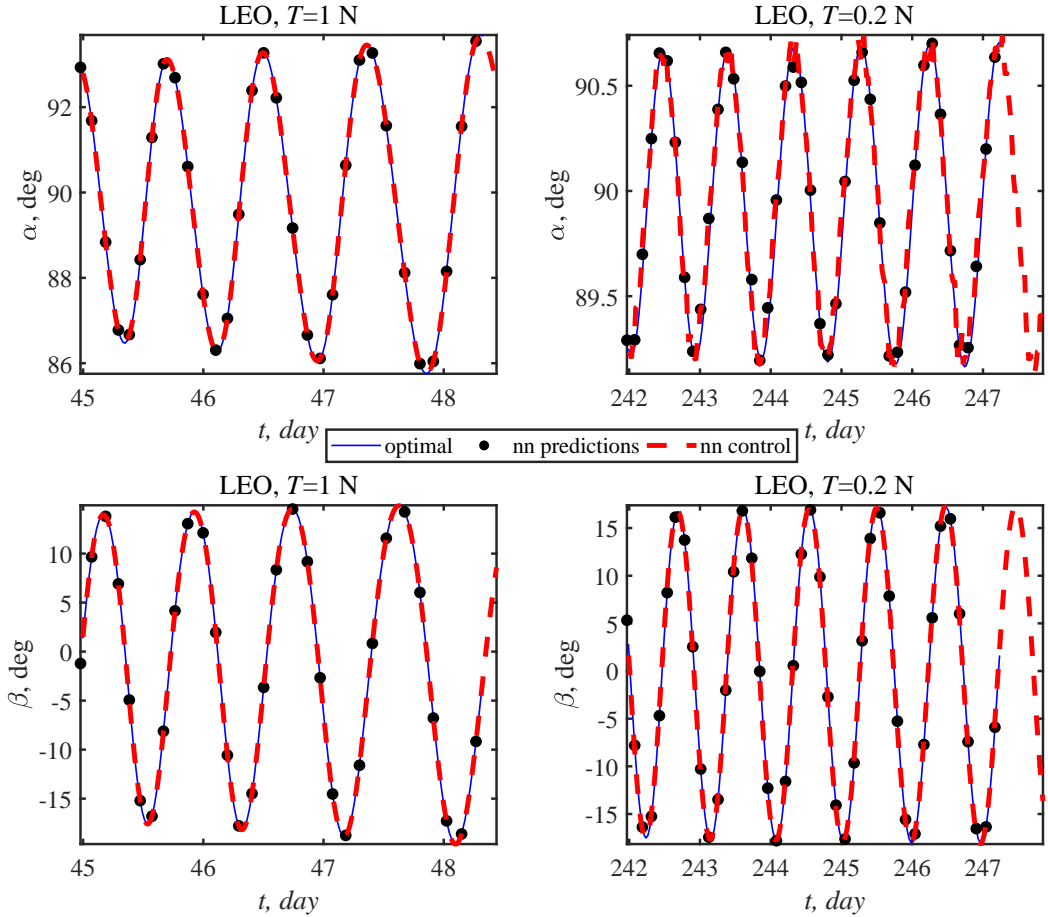}
\caption{Comparison of optimal control, NN predictions, and NN control in LEO cases}
\label{fig:LeoCompare}
\end{figure}

\begin{figure}[h]
\centering
\includegraphics[width=1.\textwidth]{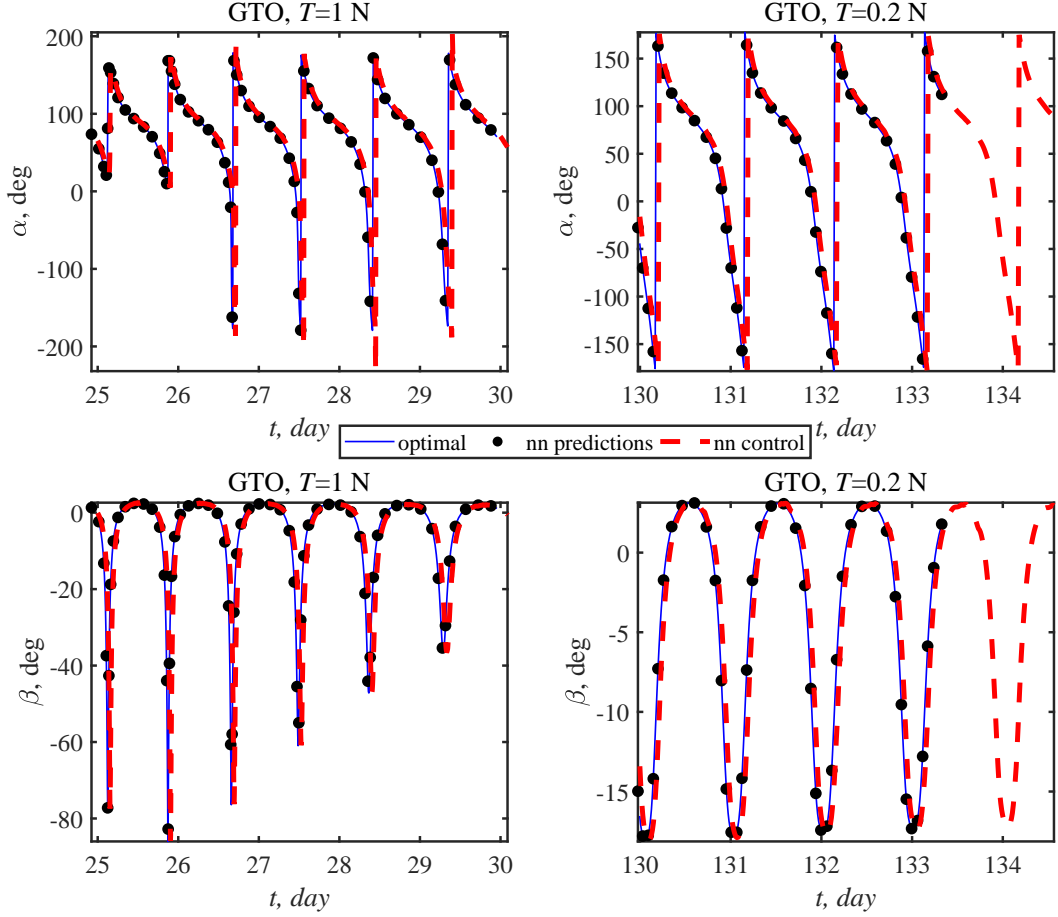}
\caption{Comparison of optimal control, NN predictions, and NN control in GTO cases}
\label{fig:GtoCompare}
\end{figure}

NN-driven trajectories from nominal departure orbits are computed, and control angles in the last several days of transfers are shown in Fig.\ref{fig:LeoCompare} and Fig.\ref{fig:GtoCompare}. The optimal control, NN predictions, and NN control are compared. The NN predictions use the states in optimal control as inputs. It can be seen that NN predictions are accurate and NN control coincides well with optimal control. In some cases the NN control continues after the optimal control stops, this is because the NN control stops at the final states that are not exactly the same as GEO. However, as the optimal mechanism has been learned, the NN can still perform the control similar to the optimal control.

\subsection{Generalization}

The generalization capability is to evaluate NN performance using data that is not in training datasets. To test the generalization capability of NNs, four test types are considered:

$Inside$ A:
A is the initialization area of departure orbits in datasets. Departure orbits are randomly created inside A using Eq.\eqref{starting}.

$Outside$ A:
Departure orbits are randomly created outside A using c$\in$[1, 1.2] in Eq.\eqref{starting}.

$Uncertainty$:
Thruster direction uncertainty is considered. The one-sigma thrust direction uncertainty is set to be 0.5 deg.

$Perturbation$:
Unknown perturbation is considered. The magnitude of the perturbation force is set to be 0.01 of the gravity, and the direction of the perturbation force is random. Perturbations such as J$_2$, atmosphere, third body, and solar pressure are not considered because they are regarded as known perturbation, and they can be included in the dataset generation process.

For each test type, one hundred random NN-driven trajectories are computed. The average of $\Delta v$ distance to GEO, final states, and optimality of NN-driven trajectories are given in Table \ref{tab_generalization}. Results show that NN-driven trajectories can be very close to GEO in all cases, especially in case "Inside A". Even though the final states are not exactly the same with GEO, NNs are considered to be capable to accomplish the orbit raising task regarding that these are long-duration many-revolution transfers and the $\Delta v$ distances are quite small. In terms of optimality, the NN-driven control is extremely close to the optimal control, which means neural networks have learned the optimal mechanism.

\begin{table}[htbp]
	\fontsize{10}{10}\selectfont
	\caption{Generalization capability of NNs}
	\label{tab_generalization}
        \centering 
    \begin{tabular}{l r r r r r r} 
        \hline\hline
        &Test type& $\Delta$, m/s & $a_f$, km & $e_f$ & $i_f$, deg & MAPE$_{opt}$\\
        \hline
        LEO, T=1 N
        &Inside A&0.4509& 42164.231186&0.000206&0.001428&2.834408e-05\%\\
        &Outside A&0.8709& 42161.778090&0.000023&0.002729&2.374521e-06\%\\
        &Perturbation&1.0482& 42171.089540&0.000157&0.002987&/\\
        &Uncertainty&2.5811& 42103.798130&0.000154&0.007023&/\\
        \hline
        LEO, T=0.2 N
        &Inside A&0.6676& 42163.143722&0.000026&0.002627&1.296593e-03\%\\
        &Outside A&1.9543& 42147.614946&0.000183&0.004002&2.056340e-03\%\\
        &Perturbation&1.4110& 42152.087963&0.000464&0.001741&/\\
        &Uncertainty&2.5920& 42095.769201&0.000028&0.005368&/\\
        \hline
        GTO, T=1 N
        &Inside A&3.9861& 42201.270567&0.000805&0.003591&3.438311e-02\%\\
        &Outside A&4.8633& 42211.150492&0.001573&0.004172&2.029110e-02\%\\
        &Perturbation&4.1550& 42216.795524&0.001118&0.002458&/\\
        &Uncertainty&6.3300& 42227.793957&0.002497&0.006605&/\\
        \hline
        GTO, T=0.2 N
        &Inside A&1.3677& 42153.839955&0.000352&0.003452&1.536073e-02\%\\
        &Outside A&2.9879& 42164.097556&0.001445&0.009086&2.070352e-03\%\\
        &Perturbation&4.0165& 42167.825382&0.001917&0.008089&/\\
        &Uncertainty&3.1512& 42166.491379&0.001488&0.008852&/\\
        \hline\hline
    \end{tabular}
\end{table}

\section{Conclusion}
In this Note, neural networks are applied for autonomous time-optimal many-revolution low-thrust orbit raising problems. 
To generate datasets, the low-thrust optimization problem is solved using the indirect method with modified equinoctial elements and thrust homotopy.
Random search is utilized to search for parameters of neural network models, and the evaluation approach of neural network performance is defined.
Simulation results show that neural networks learn the optimal mechanism well and are capable to accomplish the autonomous orbit raising task. This Note shows that besides short-duration optimal control problems that have already been studied a lot, neural networks can still have excellent performance in very-long duration optimal control problems.

\section*{Funding Sources}
Part of this research is supported by the Chinese National Natural Science Fund for Distinguished Young Scientists of China (No. 11525208).

\section*{Acknowledgments}
The first author would like to acknowledge the financial support provided by the China Scholarship Council.

\bibliography{ref}

\begin{thebibliography}{28}
\newcommand{\enquote}[1]{``#1''}
\providecommand{\natexlab}[1]{#1}
\providecommand{\url}[1]{\texttt{#1}}
\providecommand{\urlprefix}{URL }
\expandafter\ifx\csname urlstyle\endcsname\relax
  \providecommand{\doi}[1]{doi:\discretionary{}{}{}#1}\else
  \providecommand{\doi}{doi:\discretionary{}{}{}\begingroup
  \urlstyle{rm}\Url}\fi

\bibitem[{Topputo and Ceccherini(2019)}]{topputo2019catalogue}
Topputo, F., and Ceccherini, S., \enquote{A Catalogue of Parametric
  Time-Optimal Transfers for All-Electric GEO Satellites,} \emph{Modeling and
  Optimization in Space Engineering}, Springer, 2019, pp. 459--478.

\bibitem[{Ceccherini and Topputo(2018)}]{ceccherini2018system}
Ceccherini, S., and Topputo, F., \enquote{System-Trajectory Optimization of
  Hybrid Transfers to the Geostationary Orbit,} \emph{2018 Space Flight
  Mechanics Meeting}, Kissimmee, FL, 2018, p. 723.

\bibitem[{Izzo et~al.(2019{\natexlab{a}})Izzo, M{\"a}rtens, and
  Pan}]{izzo2018survey}
Izzo, D., M{\"a}rtens, M., and Pan, B., \enquote{A survey on artificial
  intelligence trends in spacecraft guidance dynamics and control,}
  \emph{Astrodynamics}, 2019{\natexlab{a}}, pp. 1--13.

\bibitem[{S{\'a}nchez-S{\'a}nchez and Izzo(2018)}]{sanchez2018real}
S{\'a}nchez-S{\'a}nchez, C., and Izzo, D., \enquote{Real-time optimal control
  via Deep Neural Networks: study on landing problems,} \emph{Journal of
  Guidance, Control, and Dynamics}, Vol.~41, No.~5, 2018, pp. 1122--1135.

\bibitem[{Cheng et~al.(2018)Cheng, Wang, Jiang, and Zhou}]{cheng2018real}
Cheng, L., Wang, Z., Jiang, F., and Zhou, C., \enquote{Real-Time Optimal
  Control for Spacecraft Orbit Transfer via Multi-Scale Deep Neural Networks,}
  \emph{IEEE Transactions on Aerospace and Electronic Systems}, 2018.

\bibitem[{Izzo et~al.(2019{\natexlab{b}})Izzo, Sprague, and
  Tailor}]{izzo2019machine}
Izzo, D., Sprague, C.~I., and Tailor, D.~V., \enquote{Machine learning and
  evolutionary techniques in interplanetary trajectory design,} \emph{Modeling
  and Optimization in Space Engineering}, Springer, 2019{\natexlab{b}}, pp.
  191--210.

\bibitem[{Furfaro et~al.(2018)Furfaro, Bloise, Orlandelli, Di~Lizia, Topputo,
  and Richard}]{furfaro2018deep}
Furfaro, R., Bloise, I., Orlandelli, M., Di~Lizia, P., Topputo, F., and
  Richard, L., \enquote{Deep Learning for Autonomous Lunar Landing,} \emph{2018
  AAS/AIAA Astrodynamics Specialist Conference}, 2018, pp. 1--22.

\bibitem[{Edelbaum(1961)}]{edelbaum1961propulsion}
Edelbaum, T.~N., \enquote{Propulsion requirements for controllable satellites,}
  \emph{ARS Journal}, Vol.~31, No.~8, 1961, pp. 1079--1089.

\bibitem[{Casalino and Colasurdo(2007)}]{casalino2007improved}
Casalino, L., and Colasurdo, G., \enquote{Improved Edelbaum's approach to
  optimize low earth/geostationary orbits low-thrust transfers,} \emph{Journal
  of guidance, control, and dynamics}, Vol.~30, No.~5, 2007, pp. 1504--1511.

\bibitem[{Kluever(2011)}]{kluever2011using}
Kluever, C.~A., \enquote{Using edelbaum's method to compute Low-Thrust
  transfers with earth-shadow eclipses,} \emph{Journal of Guidance, Control,
  and Dynamics}, Vol.~34, No.~1, 2011, pp. 300--303.

\bibitem[{Gao(2007)}]{gao2007near}
Gao, Y., \enquote{Near-optimal very low-thrust earth-orbit transfers and
  guidance schemes,} \emph{Journal of guidance, control, and dynamics},
  Vol.~30, No.~2, 2007, pp. 529--539.

\bibitem[{Guelman and Shiryaev(2016)}]{guelman2016closed}
Guelman, M.~M., and Shiryaev, A., \enquote{Closed-Loop Orbit Transfer Using
  Solar Electric Propulsion,} \emph{Journal of Guidance, Control, and
  Dynamics}, 2016, pp. 2563--2569.

\bibitem[{Zhang et~al.(2015)Zhang, Topputo, Bernelli-Zazzera, and
  Zhao}]{zhang2015low}
Zhang, C., Topputo, F., Bernelli-Zazzera, F., and Zhao, Y.-S.,
  \enquote{Low-thrust minimum-fuel optimization in the circular restricted
  three-body problem,} \emph{Journal of Guidance, Control, and Dynamics},
  Vol.~38, No.~8, 2015, pp. 1501--1510.

\bibitem[{Jiang et~al.(2017)Jiang, Tang, and Li}]{jiang2017improving}
Jiang, F., Tang, G., and Li, J., \enquote{Improving low-thrust trajectory
  optimization by adjoint estimation with shape-based path,} \emph{Journal of
  Guidance, Control, and Dynamics}, Vol.~40, No.~12, 2017, pp. 3282--3289.

\bibitem[{Chen et~al.(2018)Chen, Li, and Baoyin}]{chen2018multi}
Chen, S., Li, H., and Baoyin, H., \enquote{Multi-rendezvous low-thrust
  trajectory optimization using costate transforming and homotopic approach,}
  \emph{Astrophysics and Space Science}, Vol. 363, No.~6, 2018, p. 128.

\bibitem[{Yang et~al.(2018)Yang, Li, and Bai}]{yang2018fast}
Yang, H., Li, S., and Bai, X., \enquote{Fast homotopy method for asteroid
  landing trajectory optimization using approximate initial costates,}
  \emph{Journal of Guidance, Control, and Dynamics}, Vol.~42, No.~3, 2018, pp.
  585--597.

\bibitem[{Jiang et~al.(2012)Jiang, Baoyin, and Li}]{jiang2012practical}
Jiang, F., Baoyin, H., and Li, J., \enquote{Practical techniques for low-thrust
  trajectory optimization with homotopic approach,} \emph{Journal of Guidance,
  Control, and Dynamics}, Vol.~35, No.~1, 2012, pp. 245--258.

\bibitem[{Chi et~al.(2018)Chi, Li, Jiang, and Li}]{chi2018power}
Chi, Z., Li, H., Jiang, F., and Li, J., \enquote{Power-limited low-thrust
  trajectory optimization with operation point detection,} \emph{Astrophysics
  and Space Science}, Vol. 363, No.~6, 2018, p. 122.

\bibitem[{Gonzalo et~al.(2017{\natexlab{a}})Gonzalo, Bombardelli, and
  Topputo}]{gonzalo2017unified}
Gonzalo, J.~L., Bombardelli, C., and Topputo, F., \enquote{Unified Formulation
  for Element-Based Indirect Trajectory Optimization,} \emph{26th International
  Symposium on Space Flight Dynamics}, 2017{\natexlab{a}}, pp. 1--6.

\bibitem[{Gonzalo et~al.(2017{\natexlab{b}})Gonzalo, Topputo, and
  Armellin}]{gonzalo2017indirect}
Gonzalo, J.~L., Topputo, F., and Armellin, R., \enquote{Indirect Optimization
  of End-of-Life Disposal for Galileo Constellation Using Low Thrust
  Propulsion,} \emph{Proc. 26th International Symposium on Space Flight
  Dynamics}, 2017{\natexlab{b}}.

\bibitem[{Li et~al.(2017)Li, Chen, and Baoyin}]{li2017j2}
Li, H., Chen, S., and Baoyin, H., \enquote{J2-perturbed multitarget rendezvous
  optimization with low thrust,} \emph{Journal of Guidance, Control, and
  Dynamics}, Vol.~41, No.~3, 2017, pp. 802--808.

\bibitem[{Bergstra and Bengio(2012)}]{bergstra2012random}
Bergstra, J., and Bengio, Y., \enquote{Random search for hyper-parameter
  optimization,} \emph{Journal of Machine Learning Research}, Vol.~13, No. Feb,
  2012, pp. 281--305.

\bibitem[{Li et~al.(2019)Li, Chen, Izzo, and Baoyin}]{li2019approximators}
Li, H., Chen, S., Izzo, D., and Baoyin, H., \enquote{Deep Networks as
  Approximators of Optimal Transfers Solutions in Multitarget Missions,}
  \emph{arXiv preprint arXiv:1902.00250}, 2019.

\bibitem[{Caillau and Daoud(2012)}]{caillau2012minimum}
Caillau, J.-B., and Daoud, B., \enquote{Minimum time control of the restricted
  three-body problem,} \emph{SIAM Journal on Control and Optimization},
  Vol.~50, No.~6, 2012, pp. 3178--3202.

\bibitem[{Sutskever et~al.(2013)Sutskever, Martens, Dahl, and
  Hinton}]{sutskever2013importance}
Sutskever, I., Martens, J., Dahl, G., and Hinton, G., \enquote{On the
  importance of initialization and momentum in deep learning,}
  \emph{International conference on machine learning}, 2013, pp. 1139--1147.

\bibitem[{Kingma and Ba(2014)}]{kingma2014adam}
Kingma, D.~P., and Ba, J., \enquote{Adam: A method for stochastic
  optimization,} \emph{arXiv preprint arXiv:1412.6980}, 2014.

\bibitem[{Bengio(2012)}]{bengio2012practical}
Bengio, Y., \enquote{Practical recommendations for gradient-based training of
  deep architectures,} \emph{Neural networks: Tricks of the trade}, Springer,
  2012, pp. 437--478.

\bibitem[{Topputo and Bernelli-Zazzera(2011)}]{topputo2011optimal}
Topputo, F., and Bernelli-Zazzera, F., \enquote{Optimal low-thrust
  stationkeeping of geostationary satellites,} \emph{Proceedings of
  Confederation of European Aerospace Societies Conference}, Confine Edizioni
  Monghidoro, Bologna, 2011, pp. 1917--1925.

\end{thebibliography}

\end{document}